%

\magnification=\magstep1
\input amstex
\documentstyle{amsppt}
\NoBlackBoxes
\define\moc{{Ult}_A^{\Cal M} \Cal M}
\topmatter
\title
Splitting Number and the Core Model
\endtitle
\author
Jind\v rich Zapletal
\endauthor
\thanks
The author would like to thank Prof. T. J. Jech for infinite discussions on the 
subject and Prof. W. J. Mitchell for his interest.\endthanks
\affil
The Pennsylvania State University
\endaffil
\address
Department of Mathematics,
The Pennsylvania State University,
University Park, PA 16802
\endaddress
\email
zapletal\@math.psu.edu
\endemail
\subjclass
\endsubjclass
\abstract
We provide a lower bound for the consistency strength 
of the hypothesis proposed by S. Kamo: 
$\exists \kappa >\aleph _0$ $s(\kappa )\geq \kappa ^{++}.$
\endabstract
\endtopmatter
\document
\subhead
0.  Introduction
\endsubhead

In \cite {3} ,\cite {10} cardinal invariants on $\omega $ 
are generalized for uncountable regular cardinals and it 
is shown that some of their properties still hold true for 
their uncountable counterparts.However,the role of the splitting number
 $s(\kappa )$ (for the definition see below) changes significantly.
As Suzuki \cite {10} observed,existence of $\kappa >\aleph _0$ 
regular such that $s(\kappa )\geq \kappa ^{++} $ implies existence
 of inner models with measurable cardinals and later Kamo \cite {5}
 proved that $\exists \kappa >\aleph _0$ $s(\kappa )\geq \kappa ^{++}$
 is consistent provided there is $2^\kappa$-supercompact cardinal
 $\kappa$.Here we give a lower bound for the consistency strength
 of the statement $\exists \kappa >\aleph _0$ $s(\kappa )\geq \kappa ^{++}.$

\proclaim
{Theorem} Cons ( $\exists \kappa >\aleph _0 $ $s(\kappa )\geq \kappa ^{++}$
 ) $\to$ Cons ( $\exists \alpha $ $o(\alpha )\geq \alpha ^{++}$ ).
\endproclaim
\subhead
1.   Preliminaries
\endsubhead

\proclaim
{Definition 1} If $\kappa$ is a regular cardinal,we define the 
{\it splitting number} of $\kappa$ ,

$$s(\kappa )=min\{|\Cal S |:\Cal S \subset \Cal P(\kappa ) \text { 
\& } \forall a\in \kappa ^\kappa \text { } \exists b\in \Cal S \text {
 } |a \cap b|=|a \setminus b|=\kappa\}$$

Such $\Cal S$'s will be called {\it splitting families}.If 
$\Cal S \subset \Cal P(\kappa )$ is not splitting then there is 
a counterexample $A \in \kappa ^\kappa .$We say that $A$ {\it cuts} $\Cal S.$
\endproclaim

Now suppose $s(\kappa )\geq \kappa ^{++}$ for some $\kappa >\aleph _0$ 
regular and $M$ is an inner model satisfying GCH.Then we observe that
 $\Cal P(\kappa )\cap M$ is not a splitting family and thus there is 
$A$ cutting $\Cal P(\kappa )\cap M.$ $F=\{b\in \Cal P(\kappa ) 
\cap M:|A\cap b|=\kappa \}$ is a weakly $\kappa$-complete $M$-ultrafilter:
 if $\langle b_\alpha :\alpha <\beta <\kappa \rangle$ is a subset of $F$ 
we have $\forall \alpha <\beta $ $|A \setminus b_\alpha |<\kappa$ and so 
$|A\setminus \bigcup _{\alpha <\beta } b_\alpha |<\kappa ,$ 
$|A \cap \bigcap _{\alpha <\beta} b_\alpha |=\kappa$,in particular
  $\bigcap _{\alpha <\beta} b_\alpha \neq 0$.Hence we have 
a wellfounded ultrapower $N$ of $M$ and an elementary embedding
 $j:M\to N,$ $crit(j)=\kappa$.

These facts justify the following definitions:
\proclaim
{Definition 2} Let $\Cal M =\langle M,\in ,\lessdot ,\dots \rangle$ be 
either an inner model of set theory where $\lessdot$ is a wellordering
 of the universe,a class in $M,$or a set model of a fraction of ZFC,
not necessarily transitive,$\kappa \subset M,$where $\in$ is the "real"
 $\in$ relation and $\lessdot$ is a wellordering of $M.$For $S\subset M,$
 $S\in M$ we define

$$Def(\Cal M ,S)$$

to be the set of all subsets of $S$ first-order definable over $\Cal M$.
Obviously in the inner model case without further special predicates,
this reduces to $\Cal P(S)\cap M.$If $Def(\Cal M ,\kappa )$ is not 
a splitting family and $A$ cuts it,we define

$$F_A^\Cal M =\{a\in Def(\Cal M ,\kappa ):|a\cap A|=\kappa\}$$

As observed above,this is a weakly $\kappa$-complete filter (
or rather centered system).In this case we define

$$\moc$$

to be the ultrapower of $\Cal M$ by first-order definable over 
$\Cal M$ functions from $\kappa$ to $M$ factorized by $F_A^\Cal M.$
 \L o\' s's theorem for $\moc$ goes through due to the built-in 
wellordering and we get an elementary

$$j_A^{\Cal M }:\Cal M \to \moc$$

\endproclaim

Before we proceed to the proof of the Theorem,one easy lemma.

\proclaim
{Lemma 3} Suppose $\Cal M$ is as in the Definition 2, $A$ cuts 
$Def(\Cal M ,\kappa ),$ $S\in M,$ $S\subset M$ and $f$ is a first-order
 definable over $\Cal M$ function,$f:\kappa \to S$ , $\forall s\in S$
 $f^{-1} s\notin F_A^\Cal M.$Then

$$ F_{f^{\prime \prime}A,S}^\Cal M=\{ a\in \Cal P (S) \cap 
Def(\Cal M,S):|a\cap f^{\prime \prime}A|=\kappa \}$$

is again a weakly $\kappa$-complete $Def(\Cal M,S)$-ultrafilter.
\endproclaim
\demo
{Proof} Suppose $a\in Def(\Cal M,S)$ so $b=f^{-1}a\in Def (\Cal M,\kappa ).
$If $|a\cap f^{\prime \prime}A|<\kappa$ then the small preimages of 
singletons in $S$ guarantee $|b\cap A|<\kappa .$We get $|a\cap f^{
\prime \prime}A|=\kappa$ iff $|b\cap A|=\kappa$ iff $|A\setminus b|<\kappa$ 
iff $|f^{\prime \prime}A\setminus a|<\kappa$ and we are through.
\enddemo

\proclaim
{Corollary 4} If $\Cal M$ is as in the Definition 2 and 
$Def(\Cal M ,\kappa )$ is not a splitting family,it is possible to choose
 $A$ cutting $Def(\Cal M ,\kappa )$ such that $F_A^\Cal M$ is normal 
w.r.t. functions first-order definable over $\Cal M.$
\endproclaim
\demo
{Proof} Choose $A^\prime$ arbitrary cutting $Def(\Cal M,\kappa )$ 
and consider ${Ult}_{A^\prime}^{\Cal M} \Cal M$ and $f$ first-order
 definable over $\Cal M$ representing $\kappa .$W.l.o.g.
$f:\kappa \to \kappa .$We claim that $A=f^{\prime \prime}A^\prime$
 is just what we need.$f$ does not represent any $\alpha <\kappa ,$
and so $\forall \alpha <\kappa $ $|f^{-1}\alpha \cap A|<\kappa .$
Lemma 3 applies and $A=f^{\prime \prime}A^\prime$ cuts 
$Def (\Cal M,\kappa ).$Suppose $g$ is a $\Cal M$-definable regressive
 ordinal function,$|dom(g)\cap A|=\kappa .$Then $g\circ f$ represents
 some $\alpha <\kappa$ in ${Ult}_{A^\prime}^{\Cal M} \Cal M.$It follows
 that $|(g\circ f)^{-1}(\alpha )\cap A^\prime |=\kappa ,$by small 
preimages of singletons in $\kappa$ $|g^{-1}(\alpha )\cap f^{
\prime \prime}A^\prime |=\kappa $ and $\alpha$ is the wanted 
stabilizing value for $g.$
\enddemo

\subhead
2.  Proof of the Theorem
\endsubhead

From now on,we fix $\kappa >\aleph _0$ regular with 
$s(\kappa )\geq \kappa ^{++}$ and suppose that there are no inner models
 of $\exists \alpha$ $o(\alpha )\geq \alpha ^{++}$ .Toward the end
 of the paper,we arrive at contradiction from these assumptions,
proving the Theorem.
Recall the basic properties of Mitchell's core model $\bold K$ 
for a sequence of measures \cite {8}:
\roster
\item [1] $\bold K \models $ "ZFC,GCH,measures are wellordered,
there is a definable wellordering of the universe"
\endroster
and if there are no inner models of $\exists \alpha $ 
$o(\alpha )\geq \alpha ^{++}$ (our assumption) then
\roster
\item [2] every $j:\bold K \to M$ is an iteration of measures in $\bold K,$
\item [3] if $F$ is a $\bold K$-normal $\bold K$-ultrafilter such that
 $Ult_F \bold K$ is wellfounded then $F\in \bold K.$
\endroster

It follows from our assumptions that $\bold K \models "\kappa$ 
is measurable" . $|\Cal P (\kappa )\cap \bold K |=(\kappa ^+ )^{
\bold K} <\kappa ^{++}$ and we have $A$ cutting 
$\Cal P (\kappa )\cap \bold K$ .By Lemma 3, $A$ can be chosen such that
 $F_A^{\bold K}$ is $\bold K$-normal and hence belongs to $\bold K$ .

Still one more preparatory lemma.

\proclaim
{Lemma 5}Suppose $M$,$N$ are inner models and $j:M\to N$ is an iteration
 $\langle M_\alpha :\alpha \leq \gamma ,j_{\alpha ,\beta }:M_\alpha 
\to M_\beta :\alpha <\beta \leq \gamma \rangle$ where $M=M_0$,$N=M_\gamma$,
$j_{\alpha ,\beta}$ is a commutative system,$M_{\alpha +1}$ is 
an ultrapower of $M_\alpha$ by $M_\alpha$-measure $U_\alpha \in M_\alpha$ 
and at limit steps we take direct limits.Then if $\kappa =crit(j)$ 
we have $V_{\kappa +1}^M =V_{\kappa +1}^N$.
\endproclaim

\demo
{Proof} We prove by induction on $\alpha \leq \gamma$ that
 $V_{\kappa +1}^M =V_{\kappa +1}^{M_\alpha }.$The succesor step is 
trivial since we assume that $U_\alpha \in M_\alpha$ and
 $\kappa =min\{crit(j_{\alpha ,\alpha +1})\} .$Suppose now that 
$\alpha \leq \gamma$ is limit and $\forall \alpha^\prime <\alpha$
 $V_{\kappa +1}^M =V_{\kappa +1}^{M_{\alpha ^\prime}}$.By wellfoundedness 
of $M_\alpha ,$there were only finitely many ultrapowers taken at 
$\kappa$ and we have $\beta <\alpha$ such that 
$\forall \beta <\alpha^\prime <\alpha$ 
$crit(j_{\alpha^\prime ,\alpha^\prime +1})>\kappa.$Thus,
$V_{\kappa +1}^{M_\alpha}=dirlim_{\beta <\alpha^\prime <\alpha}
 V_{\kappa +1}^{M_{\alpha^\prime}} =V_{\kappa +1}^M.$
\enddemo

Now we proceed directly to the contradiction proving the Theorem.
The strategy is to find $\kappa ^{++}$ different normal $\bold K$-measures 
on $\kappa$.Our measures will be of the form $F_A^{\bold K}$ for some $A$.
The problem is,how should we choose $\langle A_\alpha :\alpha <\kappa ^{++} 
\rangle$ so that it is guaranteed $\alpha <\beta <\kappa ^{++} 
\to F_{A_\alpha }^{\bold K} \neq F_{A_\beta }^{\bold K}$?

Choose $\lambda \gg \kappa$ strongly limit,$\lessdot$ a wellordering of
 $H_\lambda$.Pick $A_0\subset \kappa$ such that $A_0$ cuts 
$\Cal P (\kappa )\cap \bold K$ and $F_{A_0}^{\bold K}$ 
is $\bold K$-normal.Further choose

$$\Cal M =\langle M,\in ,\lessdot ,M\cap \bold K \rangle \prec 
\langle H_\lambda ,\in ,\lessdot ,H_\lambda \cap \bold K \rangle$$

arbitrary with $A_0 \in M$,$V_{\kappa +1}^{\bold K} \subset M$ and
 $|M|\leq \kappa ^+$.This is possible since $\bold K \models "\kappa$
 is measurable" and so $|V_{\kappa +1}^{\bold K}|=(\kappa ^+)^{\bold K} 
\leq \kappa ^+$.It follows that $|Def (\Cal M ,\kappa )|\leq \kappa ^+$
 and $Def (\Cal M ,\kappa )$ is not a splitting family.Choose $A_1$ 
cutting it such that $F_{A_1}^{\Cal M}$ is normal w.r.t. functions 
first-order definable over $\Cal M$.We get

$$j_{A_1}^{\Cal M}:\Cal M \to Ult_{A_1}^{\Cal M} \Cal M 
\simeq \langle N,\in ,\lessdot ^*,N_0 \rangle$$

where N is the transitive collapse of the universe of  
$Ult_{A_1}^{\Cal M} \Cal M$.$crit (j_{A_1}^{\Cal M})=\kappa$ 
and by restriction,

$$j_{A_1}^{\Cal M}:M\cap \bold K \to N_0$$

Note that $\bold K$ is transitive and so by elementarity $N\models "N_0$
 is transitive",hence $N_0$ is "really" transitive.

\proclaim
{Claim 6} $\Cal P (\kappa )\cap N_0 =\Cal P (\kappa )\cap \bold K.$
\endproclaim

\demo
{Proof} Obviously for any $a\in \Cal P (\kappa )\cap \bold K$
 $j_{A_1}^\Cal M(a)\cap \kappa =a\in N_0$.For the other inclusion,fix 
$a\in \Cal P (\kappa )\cap N_0.$Then there is $f$ first-order definable
 over $\Cal M$ such that $[f]_{F_{A_1}^{\Cal M}}=a$.W.l.o.g. we can 
suppose $f:\kappa \to V_\kappa ^\bold K$.We distinguish two cases:
 either for some $s\in \kappa ^{<\kappa}\cap \bold K$ 
$a=j_{A_1}^\Cal M (s)=s.$Then we have $a\in \bold K$ and we are done.Or,
$\forall s\in V_\kappa ^\bold K$ $|f^{-1}s\cap A_1|<\kappa .$
In this case,define 

$$U=\{b\in V_{\kappa +1}^\bold K:a\in j_{A_1}^\Cal M (b)\} .\tag 1$$

Now by standard arguments and Lemma 3,

$$U=\{b\in V_{\kappa +1}^\bold K:|b\cap f^{\prime \prime}A_1|=\kappa\}
=F_{f^{\prime \prime}A_1,V_\kappa^\bold K}^\bold K\tag 2$$

From (2) and small preimages of singletons,$U$ is weakly $\kappa$-complete
 $\bold K$-ultrafilter. From (1),the following diagram commutes:

$$
\CD
V_\kappa ^\bold K    @>j_{A_1}^\Cal M>>   V_{j_{A_1}^\Cal M \kappa}^{N_0}\\
@\vert                                    @AAkA\\
V_\kappa ^\bold K    @>j_U>>              V_{j_U\kappa}^{N_1}
\endCD
$$      
where $j_U:\bold K\to N_1$ is given by $U$ (thus $N_1$ is wellfounded) and 
$k([g]_U)=(j_{A_1}^\Cal M g)(a)$ for any $g\in \bold K$,
$g:\kappa \to V_\kappa ^\bold K.$Note that such $g$ belongs to
 $V_{\kappa +1}^\bold K$ and the definition of $k$ makes sense.
$crit(k)\geq \kappa ,$ $k[id]_U=a$ and $[id]_U=a.$Thus 
$a\in V_{\kappa +1}^{N_1}.$But $j_U$ is an iteration of measures in 
$\bold K$ with critical point $=\kappa ,$Lemma 5 applies and we have
 $V_{\kappa +1}^{N_1}=V_{\kappa +1}^\bold K ,$thus $a\in \bold K.${\it Fini}.
\enddemo
(W. J. Mitchell pointed out a different proof of this fact using 
the fine structure of $\bold K.)$Thus $F_{A_0}^\bold K$ is a weakly 
$\kappa$-complete $N_0$-normal $N_0$-ultrafilter.Since

$$\langle H_\lambda ,H_\lambda \cap \bold K \rangle \models \text {"All 
weakly $\omega$-complete $H_\lambda \cap \bold K$-normal $H_\lambda \cap 
\bold K$-ulf's are in $H_\lambda \cap \bold K$"}$$

the same sentence is modeled in $\Cal M$ about $M\cap \bold K$ and 
by elementarity $N$ says the same about $N_0.$Now notice that 
$j_{A_1}^\Cal M (A_0)\cap \kappa =A_0\in N$ and get $F_{A_0}^\bold K\in N_0.$
At this point we observe that $F_{A_1}^\bold K =F_{A_1}^\Cal M\cap \bold K$
 and go on to prove $F_{A_1}^\bold K\neq F_{A_0}^\bold K.$

\proclaim
{Claim 7} $F_{A_1}^\bold K\neq F_{A_0}^\bold K.$
\endproclaim

\demo
{Proof} Similar as for the previous Claim.Suppose
 $F_{A_1}^\bold K=F_{A_0}^\bold K.$Then $F_{A_1}^\bold K \in N_0,$
 $F_{A_1}^\bold K=[f]_{F_{A_1}^\Cal M}$ for some
 $f:\kappa \to V_\kappa ^\bold K$ first-order definable over $\Cal M.$
As before,define $U=\{ b\in V_{\kappa +1}^\bold K:|b\cap 
f^{\prime \prime}A_1|=\kappa \}=\{ b\in V_{\kappa +1}^\bold K:
F_{A_1}^\bold K\in j_{A_1}^\Cal M (b)\}.$Again,$U$ is weakly 
$\kappa$-complete $\bold K$-ultrafilter and the following diagram commutes:
$$
\CD
V_\kappa ^\bold K    @>j_{A_1}^\Cal M>>   V_{j_{A_1}^\Cal M \kappa}^{N_0}\\
@\vert                                    @AAkA\\
V_\kappa ^\bold K    @>j_U>>              V_{j_U\kappa}^{N_2}
\endCD
$$
Here $j_U:\bold K\to N_2$ is generated by $U$ and $N_2$ is therefore 
wellfounded.$k$ is defined by $k([g]_U)=(j_{A_1}^\Cal M g)(F_{A_1}^\bold K)$
 for any $g\in \bold K$,$g:\kappa \to V_\kappa ^\bold K.$ 
$k[id]_U=F_{A_1}^\bold K$ and $crit (k)>\kappa ,$since $\kappa$ 
is definable from $F_{A_1}^\bold K.$Now $j_U$ is an iteration of measures 
in $\bold K$ ,by commutativity of iteration we can suppose that 
the (finitely many) ultrapowers on $\kappa$ were taken first.Thus 
$j_U=j_U^\prime \circ j_{U_n}\circ j_{U_{n-1}}\circ \dots \circ j_{U_0}$ 
for some $n<\omega$,$\bold K$-measures $U_i$ on $\kappa ,$ 
$U_i\in \bold K$ for $i\leq n,$ $crit(j_U^\prime) >\kappa$ and
 $j_U^\prime$ is an iteration.Set 
$N_3=j_{U_n}\circ j_{U_{n-1}}\circ \dots \circ j_{U_0}(\bold K).$
Lemma 5 applies and $V_{\kappa +1}^\bold K=V_{\kappa +1}^{N_2},$
 $V_{\kappa +2}^{N_2}=V_{\kappa +2}^{N_3}.$Now we conclude that
 $[id]_U=F_{A_1}^\bold K\in N_2.$ $(\forall a\in \Cal P(\kappa )\cap \bold K$ 
$a\in [id]_U$ iff $k(a)=a\in F_{A_1}^\bold K.)$By standard arguments, 
$V_{\kappa +2}^{N_3}=V_{\kappa +2}^{Ult_{U_n}^\bold K\bold K},$
in particular $U_n\notin N_2.$However,as $crit(j_U^\prime )>\kappa ,$
 $U_n=\{a\in \Cal P(\kappa )\cap \bold K:\kappa \in j_Ua\}.$
On the other hand,as $crit(k)>\kappa ,$ 
$F_{A_1}^\bold K=\{a\in \Cal P(\kappa )\cap \bold K:\kappa \in 
j_{A_1}^\Cal Ma\}=\{a\in \Cal P(\kappa )\cap \bold K:\kappa \in j_Ua\}=U_n,$
 $U_n\in N_2,$contradiction.
\enddemo
(In fact,this gives $F_{A_1}^\bold K>F_{A_0}^\bold K$ in the Mitchell order.)
Now it is easy to construct $\kappa ^{++}$ different weakly 
$\kappa$-complete $\bold K$-normal $\bold K$-measures on $\kappa$:
 by induction on $\alpha <\kappa ^{++}$choose $\Cal M_\alpha ,$
 $\Cal M_\alpha =\langle M_\alpha ,\in ,\lessdot ,M_\alpha \cap 
\bold K \rangle \prec \langle H_\lambda ,\in ,\lessdot ,H_\lambda \cap 
\bold K \rangle$ with $|M_\alpha |\leq \kappa ^+,$
 $V_{\kappa +1}^\bold K\subset M_\alpha,$ 
$\langle A_\beta :\beta <\alpha \rangle \subset M_\alpha$ and $A_\alpha$ 
cutting $Def(\Cal M_\alpha ,\kappa )$ such that $F_{A_\alpha }^\bold K$ 
is normal w.r.t. functions first-order definable over $\Cal M_\alpha.$
Claims 6 and 7 give that $F_{A_\alpha }^\bold K,\alpha <\kappa ^{++}$
 are pairwise distinct weakly $\kappa$-complete $\bold K$-normal 
$\bold K$-measures,as such they all belong to $\bold K$ and from 
wellordering of measures in $\bold K$ it follows that 
$o^\bold K(\kappa )\geq \kappa ^{++} \geq (\kappa ^{++})^\bold K,$
the desired contradiction with our assumption.The Theorem is proven.

\subhead
{3.   Open problems}
\endsubhead
\roster
\item
There is a wide gap between $\kappa$ supercompact and $\kappa$ with 
$o(\kappa )\geq \kappa ^{++},$calling for further consistency results.
\item
Cons ($\aleph _{\omega +2}=min\{ |\Cal S|:\Cal S\subset 
\Cal P(\aleph _\omega )$ such that
 $\forall a\in \aleph _\omega ^{\aleph _\omega }$ $\exists b\in \Cal S$
 $|a\cap b|=|a\setminus b|=\aleph _\omega \}).$
\endroster

\enddocument